\documentclass[11pt]{amsart}
\usepackage{mathrsfs}
\newtheorem{theorem}{Theorem}[section]

\theoremstyle{definition}
\newtheorem{definition}[theorem]{Definition}
\theoremstyle{remark}

\numberwithin{equation}{section}

\newtheorem{question}[theorem]{Question}

\begin{document}
\vspace{0.5in}

\title[Some New Questions on Covers and Mappings]%
{Some New Questions on Point-countable Covers and Sequence-covering Mappings}

\author{Shou Lin}
\address{Shou Lin(the corresponding author): Institute of Mathematics,
 Ningde Normal University, Ningde 352100, PR China;}
\email{shoulin60@163.com}

\author{Jing Zhang}
\address{Jing Zhang:
School of Mathematics and Statistics, Minnan Normal University,
Zhangzhou 363000, PR China}
\email{zhangjing86@126.com}
\thanks{Supported in part by the NSFC (Nos. 11801254, 11201414, 11471153) and the Training Programme
Foundation for Excellent Youth Researching Talents of Minnan Normal University (No. MJ14001).}

\subjclass[2010]{Primary 54E40; Secondary 54C10; 54E35; 54D70; 54-02}

\keywords{Metric spaces; Point-countable families; Sequence-covering mappings; Hereditarily closure-preserving families; $k$-networks; $cs$-networks}

\begin{abstract} In this survey, 37 questions on point-countable covers and sequence-covering mappings are listed, in which some of these questions have been answered. These questions are mainly related to the theory of generalized metric spaces, involving point-countable covers, sequence-covering mappings, images of metric spaces and hereditarily closure-preserving families.
\end{abstract}

\maketitle

\section{Introduction}

All the time, the problems in mathematics are a powerful driving force to lead the development of mathematics \cite{HNV04}.
Since {\it Open Problems in Topology}~\cite{vMR90} and {\it Open Problems in Topology II}~\cite{Pea07} were published, general topology and related fields have obtained huge development \cite{HMS14, HuM92, HuM02}. The new version of the book {\it Point-countable Covers and Sequence-covering Mappings}~\cite{Ls15} written by the listed author will be published soon. This book is devoted to the area of some classes of generalized metric spaces, which are defined by certain covers or networks. It systemically investigates the theory of spaces with point-countable covers and sequence-covering mappings on generalized metric spaces, including the celebrated works obtained by Chinese scholars in this field, such as point-countable covers, sequences of point-finite covers, hereditarily closure-preserving covers and star-countable covers. It is worth mentioning that there are 91 open problems were listed in this book.

In the preface to the first edition of this book~\cite{Ls15}, Arhangel'ski\v{\i} stated: ``Though reading the monograph does not require much special background, the exposition in it goes far beyond the elementary level. It contains a rich collection of deep and beautiful results of highest professional level. The book not only brings the readers to the very first line of investigations in the theory of generalized metric spaces, but also contains many intriguing and important unsolved problems, some of them old and some new. $\cdots$~I would like to mention another joyful aspect of this monograph. Its appearance marks the success of a long period of development of general topology in China, it brings to the light important contributions to the mainstream of general topology made by a very creative group of Chinese mathematicians.''

Because this book is published in Chinese, it is difficult to understand Chinese topologist' contributions and related problems in this field for the readers who are not proficient in Chinese. To provide convenience for the readers who are interested to continue in-depth study, we selects 37 questions from this book, in which some of these questions have been answered.

All spaces in this paper are assumed to be $T_2$, and all mappings are continuous and onto. Let $\mathbb{S}=\{0\}\cup\{1/n: n\in\mathbb N\}$ be the usual convergent sequence. First, recall some basic definitions.

\begin{definition} Let $X$ be a topological space, and $\mathscr{P}$ a family of subsets of $X$.
\begin{enumerate}
\item $\mathscr{P}$ is called a $k$-$network$ for $X$~\cite{OM71} if whenever $K\subset V$ with $K$ compact and $V$ open in $X$, there exists a finite
subfamily $\mathscr{P}^{\prime}\subset\mathscr{P}$ such that $K\subset \cup\mathscr{P}^{\prime}\subset V$.

\item $\mathscr{P}$ is called a $cs$-$network$ for $X$~\cite{Gu71} if, given a sequence $\{x_{n}\}$ converging to a point $x\in X$ and a neighborhood $V$ of
$x$ in $X$, then $\{x\}\cup \{x_n: n\geq n_0\}\subset P\subset V$ for some $n_0\in\mathbb N$ and some $P\in \mathscr{P}$.

\item $\mathscr{P}$ is called a $cs^*$-$network$ for $X$~\cite{Gzm87} if whenever $\{x_{n}\}$ is a sequence converging to a point $x\in V$ with $V$ open in $X$, then $\{x\}\cup \{x_{n_i}: i\in \mathbb N\}\subset P\subset V$ for some subsequence $\{x_{n_i}\}$ of $\{x_n\}$ and some $P\in \mathscr{P}$.

\item $\mathscr{P}$ is called a $wcs^*$-$network$ for $X$~\cite{LsT94} if whenever $\{x_{n}\}$ is a sequence converging to a point $x\in V$ with $V$ open in $X$, then $\{x_{n_i}: i\in \mathbb N\}\subset P\subset V$ for some subsequence $\{x_{n_i}\}$ of $\{x_n\}$ and some $P\in \mathscr{P}$.
\end{enumerate}
\end{definition}

\begin{definition}
Let $\mathscr{P}=\bigcup_{x\in
X}\mathscr{P}_{x}$ be a family of subsets of a topological space $X$ satisfying as follows: (a) for each $x\in X$, $\mathscr{P}_{x}$ is a {\it network} of $x$ in $X$, i.e., ~$\mathscr{P}_{x}\subset (\mathscr{P})_{x}$
and if $x\in G$ and $G$ is open in $X$, then $P\subset G$ for some $P\in \mathscr{P}_{x}$; (b) if $U,V\in \mathscr{P}_{x}$, then $W\subset U\cap V$ for some
$W\in \mathscr{P}_{x}$.
\begin{enumerate}
\item $\mathscr{P}$ is an $sn$-$network$ for $X$~\cite{Ls96c} if each element $P$ of $\mathscr{P}_{x}$ is a sequential neighborhood of $x$ for each $x\in X$, i.e., every sequence converging to the point $x$ is eventually in $P$. The $\mathscr{P}_{x}$ is called an $sn$-network at $x$ in $X$.

\item $\mathscr{P}$ is called a $weak$-$base$ for $X$~\cite{Arh66} if for every $G\subset X$, the set $G$ is open in $X$ whenever for each $x\in G$ there exists a $P\in \mathscr{P}_{x}$ such that $P\subset G$. The $\mathscr{P}_{x}$ is called a $weak$-$base$ at $x$ in $X$.

\item A space $X$ is $snf$-(resp. $gf$-)$countable$~\cite{Arh66, Ls97a} if $X$
has an $sn$-network (resp. a $weak$-$base$) $\mathscr{P}$ such that each $\mathscr{P}_{x}$ is countable.
\end{enumerate}
\end{definition}

The following relations hold for a topological space \cite{Ls15}:
\begin{enumerate}
\item bases $\Rightarrow$ weak-bases $\Rightarrow sn$-networks $\Rightarrow cs$-networks $\Rightarrow cs^*$-networks $\Rightarrow wcs^*$-networks.

\item bases $\Rightarrow k$-networks $\Rightarrow wcs^*$-networks.
\end{enumerate}

\begin{definition} Let $X$, $Y$ be topological spaces, and $f:X\rightarrow Y$ a mapping.
\begin{enumerate}
\item $f$ is a $compact~mapping$ (resp. $s$-$mapping$) if each $f^{-1}(y)$ is compact (resp. separable) in $X$.

\item $f$ is a $boundary$-$compact~mapping$ (resp. $at~most~boundary$-$one~mapping$) if each $\partial f^{-1}(y)$ is compact (resp. at most one point) in $X$.

\item $f$ is a $compact$-$covering~mapping$~\cite{Mic66} if each compact subset of $Y$ is the image of some compact subset of $X$.

\item $f$ is a $sequence$-$covering~mapping$~\cite{Siw71} if whenever $\{y_{n}\}$ is a convergent sequence in $Y$, there exists a convergent sequence $\{x_{n}\}$ in $X$ with
 each $x_{n}\in f^{-1}(y_{n})$.

\item $f$ is a $pseudo$-$sequence$-$covering~mapping$~\cite{GMT84, IkLT02} if whenever $S$ is a convergent sequence containing its limit point in $Y$, there is a compact subset $K$ in $X$ such that $f(K)=S$.

\item $f$ is a $sequentially~quotient~mapping$~\cite{BoS76} if whenever $\{y_{n}\}$ is a convergent sequence in $Y$, there exist a subsequence $\{y_{n_i}\}$ of $\{y_n\}$ and  a convergent sequence $\{x_i\}$ in $X$ such that
 each ~$x_{i}\in f^{-1}(y_{n_{i}})$.

\item $f$  is a 1-$sequence$-$covering~mapping$~\cite{Ls96c} if for each $y\in Y$, there is an
$x\in f^{-1}(y)$ such that whenever $\{y_{n}\}$ is a sequence converging to $y$ in $Y$, there is a sequence $\{x_{n}\}$ converging to $x$ in $X$ with each $x_{n}\in f^{-1}(y_{n})$.
\end{enumerate}
\end{definition}

The following is obvious:
\begin{enumerate}
\item 1-sequence-covering mappings $\Rightarrow$ sequence-covering mappings $\Rightarrow$ pseudo-sequence-covering mappings, and sequentially quotinet mappings.

\item boundary-compact closed mappings $\Rightarrow$ compact-covering~mappings $\Rightarrow$ pseudo-sequence-covering mappings.
\end{enumerate}

Readers may refer to R. Engelking's standard book {\it General Topology} \cite{En89} and G. Gruenhage's fine paper {\it Generalized Metric Spaces} \cite{Gr84} for unstated definitions and terminology.

\vskip 0.5cm

\section{Point-countable Covers}
Spaces with a point-countable base can be considered as the most beautiful ones with a point-countable cover. Many results on spaces with a point-countable cover can be traced back to certain properties of metric spaces, in particular, spaces having a point-countable base. Every space with a point-countable base is meta-Lindel\"{o}f. Further,
\begin{enumerate}
\item Every Fr\'{e}chet space with a point-countable $cs^*$-network is meta-Lindel\"{o}f \cite{GMT84, Ls15}.

\item Every regular Fr\'{e}chet space with a point-countable $k$-network is meta-Lindel\"{o}f~\cite{GMT84}.

\item There exists a regular space with a point-countable weak-base which is not meta-Lindel\"{o}f~\cite{GMT84}.

\item There exists a first-countable space having a $\sigma$-discrete $k$-network which is not meta-Lindel\"{o}f \cite[Example 1.4.8]{Ls15}.
\end{enumerate}

Now, we introduce a concept as follows \cite[Definition 2.1.15]{Ls15}. Let $X$ be a topological space and $\mathscr{P}$ a family of subsets of $X$. $\mathscr{P}$ is called a $cs'$-$network$ for $X$ if whenever $\{x_{n}\}$ is a sequence in $X$ converging to a point $x\in U$ with $U$ open in $X$, there exist an $n\in\mathbb N$ and a $P\in\mathscr{P}$ such that $\{x, x_n\}\subset P\subset U$. Obviously, each $cs^*$-network is a $cs'$-network in a space.

\begin{question}(\cite[Question~2.1.24]{Ls15})\quad
Is every Fr\'{e}chet space with a point-countable $cs'$-network
meta-Lindel\"{o}f?
\end{question}

The answer to Question 2.1 is in the affirmative. The following result was obtained in \cite[Corollary 4.2]{LxLs18}: Every Fr\'{e}chet space with a point-countable $cs'$-network is a hereditarily meta-Lindel\"{o}f space.

\medskip

Every $k$-space with a point-countable $k$-network is preserved by closed mappings~\cite{LsT94, Sh94}.
However, spaces with a point-countable $k$-network are not necessarily preserved by closed mappings~\cite{Sa97b}. As each compact subset of a space having a point-countable $k$-network is metrizable~\cite{GMT84}, we have the following question.

\begin{question}(\cite[Question~2.3.7]{Ls15})\quad
Let $X$ and $Y$ be topological spaces, and $f:X\rightarrow Y$ a compact-covering and closed mapping.
Does the space $Y$ have a point-countable $k$-network if $X$ has a point-countable $k$-network and each compact subset of $Y$ is metrizable?
\end{question}

Let $\mathscr{P}$ be a family of subsets of a topological space $X$. $\mathscr{P}$ is called $compact$-$countable$ if every compact subset of $X$ meets at most countably many elements of $\mathscr{P}$.

\begin{question}(\cite[Question~2.3.16]{Ls15}, \cite{LsLi09})\quad
Is every $k$-space with a compact-countable $k$-network preserved by closed mappings?
\end{question}

Michael, Nagami~\cite{MiN73} posed the following famous problem: Is every quotient $s$-image of a metric space also a compact-covering quotient and $s$-image of a metric space?  However, Chen gave a negative answer to this problem as the following two examples shows:
\begin{enumerate}
\item there is a space which is a quotient compact image of a metric space, but not any compact-covering quotient and $s$-image of a metric space~\cite{Chp99};

\item there is a regular space which is a quotient $s$-image of a metric space, but not any compact-covering quotient and $s$-image of a metric space under the condition with a $\sigma'$-set \cite{Chp03}.
\end{enumerate}
In order to characterize compact-covering $s$-images of metric spaces, we recall the concept of $cfp$-networks.

A family $\mathscr{P}$ of subsets of a topological space $X$ is called a $cfp$-$network$ for $X$~\cite{YpfL99b} if for every compact subset $K$ of $X$ and a neighborhood $V$ of $K$ in $X$ there exists a finite subfamily $\mathscr{P}'\subset\mathscr{P}$ such that $\mathscr{P}'$ can be precisely refined by a finite cover of $K$ consisting of closed subsets of $K$ and $\cup\mathscr{P}'\subset V$. Obviously, every closed $k$-network is a $cfp$-network, and every $cfp$-network is a $cs^*$-network in a topological space.

Sequential spaces with a point-countable $cs^*$-network can be characterized as quotient $s$-images of metric spaces \cite{Tan87b}, and $k$-spaces with a point-countable $cfp$-network can be characterized as compact-covering quotient and $s$-images of metric spaces \cite{Bar99, YpfL99b}. Chen's example shows that a sequential space with a point-countable $cs^*$-network may not have any point-countable $cfp$-network.

\begin{question}(\cite[Question~2.5.21]{Ls15})
\begin{enumerate}
\item Does every sequential space with a compact-countable $cs^*$-network have a point-countable $cfp$-network?

\item Is there a regular space $X$ which has a point-countable $cs^*$-network but no any point-countable $cfp$-network under (ZFC)?
\end{enumerate}
\end{question}

The answer to Question 2.4(2) is in the negative. The following result was obtained in \cite[Proposition 2.4]{LfLsSa18}: There is a compact $T_2$-space $X$ with a point-countable $cs$-network such that $X$ does not have any point-countable $cfp$-network.

\medskip

A regular space is metrizable if and only if it has a $\sigma$-compact-finite base~\cite{Bo71}. Closed images of metric spaces can be characterized by regular Fr\'{e}chet spaces with a $\sigma$-compact-finite $k$-network \cite{Lc92}. It remains an open problem \cite{LcT07} whether a regular space with a $\sigma$-compact-finite weak-base has a $\sigma$-locally finite weak-base. Banakh, Bogachev and Kolesnikov~\cite{BaBoKo08} proved that a separable regular space with a $\sigma$-compact-finite $k$-network has a countable $k$-network under (CH).

\begin{question}(\cite[Question~4.1.24]{Ls15})\quad
Does every separable regular space with a $\sigma$-compact-finite weak-base have a countable weak-base under (ZFC)?
\end{question}

A compact closed mapping is called a {\it perfect mapping}~\cite{En89}. A countable-to-one perfect image of a space with a countable weak-base is not necessarily $gf$-countable.

\begin{question}(\cite[Question~2.6.24]{Ls15}, \cite{YLY10})\quad
Is every space with a point-countable weak-base preserved by finite-to-one closed mappings?
\end{question}

Michael and Nagami~\cite{MiN73} showed that each compact subset of a space with a point-countable base has a countable outer base.
Let $\mathscr{P}=\bigcup_{x\in X}\mathscr{P}_{x}$  be a weak-base for a space $X$ in Definition 1.2(2) and $A\subset X$. The family $\bigcup_{x\in A}\mathscr{P}_{x}$ is called an {\it outer weak-base} \cite{LsZha13} of the set $A$ in $X$.

\begin{question}(\cite[Question~2.7.20]{Ls15}, \cite{LsZha13})\quad
Does every compact subset of a space with a point-countable weak-base have a countable outer weak-base?
\end{question}

The answer to Question 2.7 is in the affirmative. The following result was obtained in \cite[Proposition 2.12]{LfLsSa18} and \cite[Lemma 2]{Tu14}: Let $X$ be a space with a point-countable weak-base. Then every compact subset of $X$ has a countable outer weak-base.

\medskip

Let $\mathscr{P}$ be a cover of a topological space $X$. $\mathscr{P}$ is called {\it point-regular} \cite{Ale60} if $x\in U$ and $U$ is open in $X$, then $\{P\in(\mathscr{P})_{x}: P\not\subset U\}$ is finite.
\begin{enumerate}
\item Spaces with a point-regular base can be characterized as (compact-covering) open and compact images of metric spaces \cite{Arh62}.

\item Spaces with a point-regular weak-base can be characterized as sequence-covering quotient and compact images of metric spaces \cite{IkLT02}.

\item Spaces with a point-regular $cs$-network (resp. $sn$-network) can be characterized as sequence-covering compact images of metric spaces \cite{LsY01b}.

\item Spaces with a point-regular $cs^*$-network can be characterized as pseudo-sequence-covering (and sequentially quotient)
 compact images of metric spaces \cite{AnTu11}.
\end{enumerate}

\begin{question}(\cite[Question~3.4.20]{Ls15}, \cite{LsZhu13})\quad
How to characterize spaces with a point-regular $k$-network by certain images of metric spaces?
\end{question}

Shakhmatov~\cite{Sha84} and Watson~\cite{Ws85} constructed independently the following famous example. There exists a completely regular pseudo-compact non-compact space $X$ having a point-countable base. Uspenskii~\cite{Us84} proved that pseudo-compact spaces with a $\sigma$-point-finite base is metrizable.

\begin{question}(\cite[Question~3.4.21]{Ls15})\quad
Is every pseudo-compact space with a point-regular weak-base metrizable?
\end{question}

\vskip 0.5cm
\section{Sequence-covering mappings}
A finite-to-one closed mapping defined on a metric space is not necessarily a sequence-covering mapping \cite{Siw71}.
Closed mappings on regular spaces in which each point is a $G_\delta$-set are sequentially quotient mappings \cite[Lemma 2.3.3]{Ls15}.

\begin{question}(\cite[Question~2.3.15]{Ls15}, \cite{Ls99b})\quad
Is every closed mapping on spaces in which each point is a $G_\delta$-set sequentially quotient?
\end{question}

The answer to Question 3.1 is in the negative. The following result was obtained in \cite[Proposition 2.2]{LfLsSa18}: There is a closed mapping $f: X\to\mathbb S$ which is not sequentially quotient such that $X$ is $T_2$ (non-regular) and every point of $X$ is a $G_\delta$-set.

\medskip

Siwice~\cite{Siw71} showed that every open mapping on a first-countable space is a sequence-covering mapping. In fact, every almost-open mapping on a first-countable space is a 1-sequence-covering mapping \cite{Ls00a}. Let $X$ and $Y$ be topological spaces, and $f:X\to Y$ a mapping.
$f$ is called an {\it almost-open mapping} \cite{Arh62a} if for each $y \in Y$, there is an $x\in f^{-1}(y)$ such that
whenever $U$ is a neighborhood of $x$ in $X$, then $f(U)$ is a neighborhood of $y$
in $Y$. There exists an open mapping on a Fr\'echet space such that it is not pseudo-sequence-covering \cite{YLY10}.

A space $X$ is {\it strongly Fr\'{e}chet} at a point $x\in X$ \cite{Siw71}
if for each decreasing sequence $\{A_{n}\}$ of subsets of $X$ with $x\in
\bigcap_{n\in\mathbb{N}}\overline{A_{n}}$, there is an $x_{n}\in
A_{n}$ for each $n\in\mathbb{N}$
such that the sequence $\{x_{n}\}$ converges to the point $x$. A space $X$ is strongly Fr\'{e}chet \cite{Siw71} if it is
strongly Fr\'{e}chet at each point of $X$.

\begin{question}(\cite[Question~2.6.19]{Ls15}, \cite{LsZhu13})\quad
Is each almost-open mapping on a strongly Fr\'{e}chet space a sequence-covering mapping?
\end{question}

The answer to Question 3.2 is in the negative. The following result was obtained in \cite[Theorem 2.6]{LfLsSa18}: There is an open mapping $\varphi: X\to\mathbb S$ which is not sequence-covering such that $X$ is countable and bi-sequential.

\medskip

Sequence-covering mappings are 1-sequence-covering mappings under certain conditions. For example, boundary-compact sequence-covering mappings on first-countable spaces are 1-sequence-covering mappings \cite[Theorem 3.5.3]{Ls15}, in particularly, compact sequence-covering mappings on metric spaces are
1-sequence-covering mappings \cite{LsY01a}.

A regular space $X$ is called $g$-$metrizable$ \cite{Siw74} if $X$ has a $\sigma$-locally finite weak-base.

\begin{question}(\cite[Question~3.5.8]{Ls15}, \cite{LsC16})\quad
Let $X$ and $Y$ be topological spaces, and $f:X\to Y$ a boundary-compact sequence-covering mapping.
Is $f$ a 1-sequence-covering mapping if $X$ satisfies one of the following conditions?
\begin{enumerate}
\item every compact subset of $X$ has a countable $sn$-network in $X$.

\item every compact subset of $X$ has a countable outer $sn$-network in $X$.

\item $X$ has a compact-countable $sn$-network.

\item $X$ is a $g$-metrizable space.
\end{enumerate}
\end{question}

The following relations hold in the conditions above: $(4)\Rightarrow (3)\Rightarrow (2)$.
The answers to the space $X$ satisfying the conditions (2)-(4) in Question 3.3 are in the negative. The following result was obtained in \cite[Proposition 2.14]{LfLsSa18} and \cite[Example 3.5]{LxLc19}: There is a boundary-compact, sequence-covering mapping $\varphi: X\to Y$ which is not 1-sequence-covering such that $X$ is $g$-second countable.

\medskip

Spaces with a point-countable base are preserved by countably bi-quotient $s$-mappings \cite{Fil69}. What kind of mappings can preserve spaces with a
point-countable $sn$-network (resp. weak-base)? Spaces with a point-countable $sn$-network are preserved by countable-to-one 1-sequence-covering mappings \cite{LfLs10}. Moreover, spaces with a point-countable $sn$-network can be characterized as 1-sequence-covering, at most boundary-one and $s$-images of metric spaces \cite[Theorem 2.6.12]{Ls15}.

\begin{question}(\cite[Question~2.6.21]{Ls15}, \cite{LfLs10})\quad
Is every space with a point-countable $sn$-network preserved by 1-sequence-covering, at most boundary-one and $s$-mappings?
\end{question}

The answer to Question 3.4 is in the negative. The following result was obtained in \cite[Theorem 2.8]{LfLsSa18}: There is a 1-sequence-covering, at most boundary-one and $s$-mapping $\varphi: X\to Y$ such that $X$ has a $\sigma$-point-finite weak-base, but $Y$ does not have any
point-countable $sn$-network. The result also answers the following question in the negative (\cite[Question~2.6.23]{Ls15}, \cite{Lc07}): Is every space with a point-countable weak-base preserved by quotient, at most boundary-one and $s$-mappings?

\medskip

Sequence-covering and closed images of regular spaces with a point-countable weak-base are $gf$-countable \cite{Lc05}.

\begin{question}(\cite[Question 3.5.19]{Ls15}, \cite{LfLs14})\quad
Let $f:X\to Y$  be a sequence-covering and closed mapping. Is $Y$ an $snf$-countable space if $X$ has a point-countable $sn$-network?
\end{question}

Sequence-covering and closed images of regular spaces with a $\sigma$-compact-finite weak-base also have a $\sigma$-compact-finite weak-base \cite{Lc05}.

\begin{question}(\cite[Question 4.1.29]{Ls15}, \cite{LfSrx10})\quad
Does every sequence-covering and closed image of a space with a $\sigma$-compact-finite $sn$-network have a $\sigma$-compact-finite $sn$-network?
\end{question}

The answers to Questions 3.5 and 3.6 are in the negative. The following result was obtained in \cite[Theorem 2.18]{LfLsSa18}: There is a sequence-covering closed mapping $\varphi: Y\to S_{\omega}$ such that $Y$ has a $\sigma$-compact-finite $sn$-network. The sequential fan $S_\omega$ is not $snf$-countable \cite[Lemma 2.17]{LfLsSa18}.

\medskip

Regular spaces with a point-countable base are preserved by sequence-covering and closed mappings \cite{Lc07a}.

\begin{question}(\cite[Question~3.5.18]{Ls15})\quad
Let $f:X\to Y$ be a sequence-covering and closed mapping. Does $Y$ have a point-countable weak-base if $X$ is a regular space with a point-countable weak-base?
\end{question}

We shall strengthen the forms of compact-covering and sequence-covering (resp. 1-sequence-covering) mappings. Let $f:X\rightarrow Y$ be a mapping. $f$ is called a 1-$scc$-$mapping$ \cite{LsZha13} if for each compact subset $K\subset Y$, there is a compact subset $L\subset X$ such that $f(L)=K$ and for each $y\in K$, there is a point $x\in L$ such that if whenever $\{y_{n}\}$ is a sequence converging to $y$, there exists a sequence $\{x_n\}$ converging to $x$ in $X$ with each~$x_n\in f^{-1}(y_n)$.
$f$ is called an $scc$-$mapping$ \cite{LsZha13} if for each compact subset $K\subset Y$, there is a compact subset $L\subset X$ such that $f(L)=K$ and if whenever $\{y_{n}\}$ is a sequence converging to some point in $K$, there exists a sequence $\{x_n\}$ converging to some point in $L$ with each $x_n\in f^{-1}(y_n)$.

The following are known that
\begin{enumerate}
\item Every 1-$scc$-mapping (resp. $scc$-mapping) is a 1-sequence-covering (resp. sequence-covering) and compact-covering mapping.

\item Every compact-covering open mapping on a first-countable space is a 1-$scc$-mapping \cite{LsZha13}.

\item Every $scc$-mapping on a first-countable space is a 1-$scc$-mapping \cite[Theorem 2.7.15]{Ls15}.

\item There exists a 1-sequence-covering compact-covering mapping on a metric space which is not an $scc$-mapping \cite{LsZha13}.
\end{enumerate}

\begin{question}(\cite[Question~2.7.16]{Ls15}, \cite{LsZha13})\quad
Is every $scc$-mapping on compact spaces a 1-$scc$-mapping?
\end{question}

The answer to Question 3.8 is in the negative. The following result was obtained in \cite[Theorem 2.10]{LfLsSa18}: Not every $scc$-mapping of a compact space is 1-$scc$.

\medskip

Sequence-covering mappings are pseudo-sequence-covering mappings and sequentially quotient mappings. For a topological space $X$ every sequence-covering mapping onto the space $X$ is a 1-sequence-covering mapping if and only if for each point $x\in X$, there is a sequence $\{x_n\}$ converging to $x$ in $X$ such that the set $\{x_n: n\in\mathbb N \}\cup\{x\}$ is a sequential neighborhood of $x$ in $X$ \cite{LsLiGe13}. On the other hand, every sequentially quotient mapping onto a space $X$ is a 1-sequence-covering mapping if and only if there is a nontrivial convergent sequence in $X$ \cite{LsLiGe13}.

\begin{question}(\cite[Question~3.5.27]{Ls15})\quad
Give a characterization of a space $X$ such that every pseudo-sequence-covering mapping onto the space $X$ is a 1-sequence-covering mapping.
\end{question}

The following are some answers to Question 3.9. For a space $X$, the following are equivalent \cite[Proposition 2.20]{LfLsSa18}.

(1)\, Every pseudo-sequence-covering mapping onto $X$ is 1-sequence-covering;

(2)\, Every pseudo-sequence-covering mapping onto $X$ is sequence-covering;

(3)\, Every convergent sequence of $X$ is a finite set;

(4)\, Every mapping onto $X$ is 1-sequence-covering.

\vskip 0.5cm
\section{Images of metric spaces}
In the past fifty years, many interesting works are induced by certain mappings on metric spaces. Obviously, every pseudo-sequence-covering mapping on a metric space is a sequentially quotient mapping \cite[Lemma 1.3.4]{Ls15}. Every sequentially quotient and boundary-compact mapping on a developable space or a space with a point-countable base is a pseudo-sequence-covering mapping \cite{LfLs14a}.

\begin{question}(\cite[Question~3.5.21]{Ls15})\quad
Is every sequentially quotient $s$-mapping on a metric space a pseudo-sequence-covering mapping?
\end{question}

Every $snf$-countable space can be characterized as a sequentially quotient and boundary-compact image of a metric space \cite{LfLs10}.
The sequentially quotient, at most boundary-one and $s$-image of a metric space can be characterized as a space with a point-countable $sn$-network \cite{LfLs10}.

\begin{question}(\cite[Question~2.6.22]{Ls15}, \cite{LfLs10})\quad
Let $X$ be a sequentially quotient $s$-image of a metric space. Is $X$ a sequentially quotient, boundary-compact and $s$-image of a metric space if $X$ is $snf$-countable?
\end{question}

A topological space $X$ is called {\it feebly compact} \cite{Ste77} if every locally finite family of open subsets of $X$ is finite. It is clear that countably compact spaces are feebly compact, and that feebly compact spaces are pseudo-compact. Arhangel'ski\v{\i}~\cite{Arh11} showed that if a regular feebly compact space $X$ is a pseudo-open $s$-image of a metric space, then $X$ has a point-countable base. If a countably compact space $X$ is a quotient $s$-image of a metric space, then $X$ is metrizable \cite{GMT84}. Here, a mapping $f: X\to Y$ is called {\it pseudo-open} \cite{Arh63} if $V$ is an open subset of $X$ and $f^{-1}(y)\subset
V$ for some $y\in Y$, then $f(V)$ is a neighborhood of $y$ in $Y$. Obviously, every pseudo-open mapping is a quotient mapping.

\begin{question}(\cite[Question~2.2.12]{Ls15})\quad
Let $X$ be a feebly compact space. Does $X$ have a point-countable base if $X$ is a quotient $s$-image of a metric space?
\end{question}

Liu \cite{Lc07} proved that a space is a quotient, at most boundary-one and $s$-image of a metric space if and only if it is a quotient, at most boundary-one and countable-to-one image of a metric space.

\begin{question}(\cite[Question~3.4.14]{Ls15}, \cite{LfLs10})\quad
Is every quotient compact image of a metric space a countable-to-one quotient image of a metric space?
\end{question}

While locally separable metric spaces are between separable metric spaces and metric spaces, the images of locally separable metric spaces are very different from the images of separable metric spaces, or the images of metric spaces \cite[Example 5.1.23]{Ls15}. Sequentially quotient $s$-images of metric spaces are equivalent to pseudo-sequence-covering $s$-images of metric spaces.

\begin{question}(\cite[Question~5.1.8]{Ls15})\quad
Is every sequentially quotient $s$-image of a locally separable metric space equivalent to a pseudo-sequence-covering $s$-image of a locally separable metric space?
\end{question}

Although some intrinsic characterizations about quotient $s$-images of locally separable metric spaces have been obtained, most of these characterizations are more complicated \cite[Theorems 5.1.4, 5.1.7 and 5.1.9]{Ls15}. It's still an open question: how to seek a nice intrinsic characterization of a quotient $s$-image of a locally separable metric space \cite{LcT96b, TaX96}. Some attempts have been made to solve the above question by the following questions.

\begin{question}(\cite[Question~5.1.15]{Ls15}, \cite{LsCL09})\quad
Is each closed Lindel\"of subspace of a space $X$ is separable if the space $X$ is a quotient $s$-image of a locally separable metric space?
\end{question}

Every first-countable subspace of a quotient $s$-image of a locally separable metric space is locally separable. However, there exists a space $X$ which is a quotient $s$-image of a metric space such that every first-countable subspace of the space $X$ is locally separable, but $X$ is not any quotient $s$-image of a locally separable metric space \cite{LsSa07, Sa06}.

\begin{question}(\cite[Question~5.1.16]{Ls15}, \cite{LsCL09})\quad
Let $X$ be a quotient $s$-image of a metric space. If every first-countable subspace of $X$ is locally separable and every closed Lindel\"of subspace of $X$ is separable, is $X$ a quotient $s$-image of a locally separable metric space?
\end{question}

A topological space $X$ is called a {\it cosmic space}~(resp. an {\it $\aleph_0$-space}) \cite{Mic66} if it has a countable network (resp. $k$-network). Cosmic spaces and $\aleph_0$-spaces are very nice spaces represented as certain images of metric spaces. Michael~\cite{Mic66} proved that a space is a cosmic space if and only if it is an image of separable metric space, and that a space is an $\aleph_0$-space if and only if it is a compact-covering image of separable metric space.

It is clear that every quotient $s$-image of a locally separable metric space is a sequential space with a point-countable $cs^*$-network consisting of~cosmic subspaces.

\begin{question}(\cite[Question~5.1.24]{Ls15}, \cite{LsCL09})\quad
Is every sequential space with a point-countable $cs^{*}$-network consisting of cosmic subspaces a quotient $s$-image of a locally separable metric space? \end{question}

Let $\mathscr{P}$  be a cover of a topological space $X$.  $\mathscr{P}$ is called a {\it $cs^{*}$-cover} of the space $X$ \cite{Ljj00a} if, for each convergent sequence $S$ in $X$, there exists a $P\in\mathscr{P}$ such that some subsequence of $S$ is frequently in $P$. $\mathscr{P}$ is called a {\it $k$-cover} of the space $X$ \cite{McN88} if every compact subset of $X$ is covered by some finite subfamily of $\mathscr{P}$.
Spaces with a point-countable $cs^*$-cover (resp. $k$-cover) consisting of $\aleph_0$-subspaces are pseudo-sequence-covering (resp. compact-covering) $s$-images of locally separable metric spaces.

\begin{question}(\cite[Question~5.1.25]{Ls15}, \cite{Ljj00c, LiJi02})
\begin{enumerate}
\item Is every pseudo-sequence-covering $s$-image of locally separable metric spaces a space with a point-countable $cs^*$-cover consisting of $\aleph_0$-subspaces?  It is even unknown whether every quotient $s$-image of locally separable metric spaces has a point-countable $cs^*$-network consisting of $\aleph_0$-subspaces \cite{LLD97, LcT96b}.

\item Is every compact-covering $s$-image of locally separable metric spaces a space with a point-countable $k$-cover consisting of $\aleph_0$-subspaces?
\end{enumerate}
\end{question}

Sequentially quotient compact images of locally separable metric spaces are equivalent to pseudo-sequence-covering compact images of locally separable metric spaces \cite{Ge03}. A characterization of the images has been obtained, but it is quite complex \cite{LvLi05}. Let $X$ be a topological space and $\{\mathscr{P}_n\}$ a sequence of covers of the space $X$. $\{\mathscr{P}_n\}$ is called a {\it point-star network} \cite{LsY01b} for $X$ if the family $\{\mbox{st}(x, \mathscr{P}_n): n\in\mathbb N\}$ is a network of $x$ in $X$ for each $x\in X$.

\begin{question}(\cite[Question~5.2.4]{Ls15})
\begin{enumerate}
\item Let $X$ be a space with a point-star network of point-finite $cs^*$-covers consisting of cosmic subspaces. Is $X$ a pseudo-sequence-covering compact image of a locally separable metric space \cite{LsCL09}?

\item  Let $X$ be a sequential space with a point-star network of point-finite $cs^*$-covers consisting of cosmic subspaces. Is $X$ a quotient compact image of a locally separable metric space \cite{GeLs04}?
\end{enumerate}
\end{question}

Let $(X,d)$ be a metric space and $Y$ a topological space. A mapping $f: X\rightarrow Y$ is called a {\it $\pi$-mapping} \cite{Pon60} if $d(f^{-1}(y), X\setminus f^{-1}(U))>0$ for each $y\in Y$ and each neighborhood $U$ of $y$ in $Y$. Each compact mapping is a $\pi$-mapping on a metric space.

\begin{question}(\cite[Question~5.2.11]{Ls15}, \cite{AnDu08})\quad
Is every quotient $\pi$-image of a locally separable metric space equivalent to a pseudo-sequence-covering quotient $\pi$-image of a locally separable metric space?
\end{question}

Quotient compact images of locally separable metric spaces are not necessarily quotient compact and compact-covering images of metric spaces \cite{Chp99}.
Quotient $\pi$-images of separable metric spaces are pseudo-sequence-covering quotient and compact images of separable metric spaces \cite{AnTu12}.
Quotient compact and regular images of separable metric spaces are quotient compact and compact-covering images of separable metric spaces \cite{LsY01a}.

\begin{question}(\cite[Question~3.3.24]{Ls15})\quad
Is every quotient compact image of a separable metric space a quotient compact and compact-covering image of a separable metric space?
\end{question}

\vskip 0.5cm
\section{Hereditarily closure-preserving families}
The research about hereditarily closure-preserving families is mainly originated from closed images of locally finite families in a topological space. Let $\mathscr{P}$ be a family of subsets of a topological space $X$. $\mathscr{P}$ is called a {\it hereditarily closure-preserving family} \cite{Las66} of the space $X$ if the family $\{H(P): P\in\mathscr{P}\}$ is closure-preserved for each $H(P)\subset P\in\mathscr{P}$, i.e., $$\overline{\cup\{H(P): P\in\mathscr{P}'\}}=\cup\{\overline{H(P)}: P\in\mathscr{P}'\}$$ for each $\mathscr{P}'\subset\mathscr{P}$.
$\mathscr{P}$ is called a {\it point-discrete family} \cite{LcLsLu08} or a {\it weakly hereditarily closure-preserving family} \cite{BEL75} of the space $X$ if the family $\{\{p(P)\}: P \in\mathscr{P}\}$ is closure-preserving for each $p(P)\in P\in\mathscr{P}$, i.e., the set $\{p(P): P\in\mathscr{P}\}$ is closed discrete in $X$. Obviously, every locally finite family is hereditarily closure-preserving, and every hereditarily closure-preserving family is point-discrete.

In regular spaces, the closure $\overline{\mathscr{P}}=\{\overline{P}: P\in\mathscr{P}\}$ of a hereditarily closure-preserving family $\mathscr{P}$ is still hereditarily closure-preserving \cite{Ls88a, MOS76}.
Generally, we know the closure of a hereditarily closure-preserving family is point-discrete \cite{SGG06}.

\begin{question}(\cite[Question~4.2.7]{Ls15}, \cite{Ls99b})\quad
Is the closure of every hereditarily closure-preserving family still hereditarily closure-preserving?
\end{question}

The answer to Question 5.1 is in the negative. The following result was obtained in \cite[Proposition 2.22]{LfLsSa18}: There are a $T_2$ (non-regular) space $X$ and a hereditarily closure-preserving family $\mathscr{P}$ in $X$ such that $\overline{\mathscr{P}}=\{\overline{P}: P\in\mathscr{P}\}$ is not hereditarily closure-preserving.

\medskip

About the relations between a $\sigma$-point-discrete family and a $\sigma$-compact-finite family, it is known that a space with a $\sigma$-point-discrete $sn$-network (resp. $k$-network, $wcs^*$-network, network) has a $\sigma$-compact-finite $sn$-network \cite{LfSrx10}(resp. $k$-network \cite{QzbG99}, $wcs^*$-network \cite{SrxLs13}, network \cite{Ls15}). However, there exists a space with a $\sigma$-point-discrete base (resp. weak-base) which does not have a $\sigma$-compact-finite base \cite{BEL75} (resp. weak-base \cite{LsYl04}).

\begin{question}(\cite[Question~4.1.19]{Ls15})\quad
Does every space with a $\sigma$-point-discrete $cs$-network have a $\sigma$-compact-finite $cs$-network \cite{SrxLs13} or a $\sigma$-compact-finite $cs^*$-network?
\end{question}

Sometimes, an $snf$-countable space having certain $cs$-networks may imply that it has certain $sn$-networks. For example, a space with a point-countable $sn$-network if and only if it is an $snf$-countable space with a point-countable $cs$-network \cite{LsY01c}.

\begin{question}(\cite[Question~4.1.20]{Ls15})\quad
Does every $snf$-countable space with a $\sigma$-point-discrete $cs$-network have a $\sigma$-point-discrete $sn$-network?
\end{question}

Because an $snf$-countable space with a $\sigma$-point-discrete $cs$-network has a $\sigma$-compact-finite $sn$-network \cite{LfLs14, LfSrx10}, a $gf$-countable space with a $\sigma$-point-discrete $cs$-network has a $\sigma$-compact-finite weak-base. However, there exists a $gf$-countable space with a $\sigma$-point-discrete $wcs^*$-network which is not a space having a $\sigma$-point-discrete $cs^*$-network \cite[Example~4.1.16(7)]{Ls15}.

\begin{question}(\cite[Question~4.1.22]{Ls15}, \cite{SrxLs13})\quad
Does every $gf$-countable space with a $\sigma$-point-discrete $cs^*$-network have a $\sigma$-compact-finite weak-base?
\end{question}

Liu and Ludwig \cite{LcLu05} proved each closed mapping on a regular space with a $\sigma$-point-discrete base is a compact-covering mapping. Liu, Lin and Ludwig \cite{LcLsLu08} proved that every closed mapping on a regular space with a $\sigma$-point-discrete weak-base is also a compact-covering mapping under (CH).

\begin{question}(\cite[Question~4.1.27]{Ls15})\quad
Is every closed mapping on a regular space with a $\sigma$-point-discrete $sn$-network a compact-covering mapping?
\end{question}

\begin{question}(\cite[Question~4.1.28]{Ls15}, \cite{LfSrx10})\quad
Is every space with a $\sigma$-point-discrete $sn$-network preserved by a sequence-covering closed mapping?
\end{question}

Lin, Liu and Dai \cite{LLD97} proved that a regular space $X$ with a $\sigma$-hereditarily closure-preserving $k$-network consisting of $\aleph_0$-subspaces if and only if $X$ has a $\sigma$-hereditarily closure-preserving $k$-network and each first-countable closed subspace of $X$ is locally separable.

\begin{question}(\cite[Question~4.3.11]{Ls15})\quad
Is each first-countable subset of a regular space with a $\sigma$-hereditarily closure-preserving $k$-network consisting of separable subspaces locally separable?
\end{question}

\vskip 0.9cm
\bibliographystyle{plain}

\end{document}